\documentclass[12pt]{amsart}

\usepackage{amsmath,amssymb,amsfonts}
\usepackage{kpfonts}
\usepackage{mathtools}
\usepackage{bbm}
\usepackage{hyperref}
\usepackage[foot]{amsaddr}

\newtheorem{theorem}{Theorem}

\newtheorem{lemma}{Lemma}

\newtheorem{prop}{Proposition}

\def\endpf{\hfill $\Box$ \vskip0.5cm}
\def \proof{\vspace{2mm}\noindent{\it Proof.\ }}

\begin{document}

\title[{Records  in the Infinite Occupancy Scheme}]{Records  in the Infinite Occupancy Scheme}

\author{Zakaria Derbazi$^1$, Alexander Gnedin$^1$}
\address{$^1$School of Mathematical Sciences, Queen Mary University of London, London E14NS, United Kingdom}

\author{Alexander Marynych$^2$}
\address{$^2$Faculty of Computer Science and Cybernetics, Taras Shev\-chen\-ko National University of Kyiv, Kyiv 01601, Ukraine}

\begin{abstract}
\noindent
We consider the classic infinite occupancy scheme, where balls are thrown  in boxes independently, with probability $p_j$ of hitting box $j$. Each time a box receives its first ball we speak of a {\it record} and, more generally, 
call an {\it $r$-record} every event when a box receives its $r$th ball. 
Assuming that the sequence $(p_j)$ is not decaying too fast, we show that after many 
balls have been thrown, the suitably scaled point  process 
of $r$-record times is approximately Poisson.  
The joint convergence of $r$-record processes is argued under a condition of regular variation.

\end{abstract}



\maketitle

\section{Introduction}
In the infinite occupancy scheme first systematically studied by Karlin~\cite{Karlin}, balls are  allocated independently  to an infinite series of boxes, with fixed probability $p_j$ of hitting the $j$th box for each ball.
There is  extensive literature on asymptotic properties of  the random partition 
associated with the allocation of a large number of balls. 
The most explored features include  the number of boxes occupied by at least one ball, and the counts of boxes occupied $r$ times.
See~\cite{ GBarbour, Ben, GHP, Hwang} for development and many references therein.
Much less attention has been devoted to the evolutionary aspects of the partition seen as a random process when balls are thrown
successively  one at a time.

Recently functional Gaussian limits were shown 
for the mentioned statistics
under Karlin's condition  of  regular variation on $(p_j)$~\cite{CZ, DW}.
However,  
the contracted time scale employed for such approximation
turns  too rough to apprehend a short term pattern of newly `discovered' boxes, that get hit for the first time.  
The latter aspect is of interest for statistical applications akin to the new species search  problem \cite{Bunge}.

Following the terminology from~\cite{Spano},  in this paper we call a
{\it record}  each occupancy event  when a new box is hit and, more generally,  call an {\it $r$-record} each event when a box is hit for the $r$th time.
The class of models considered  here has the property that the number of boxes occupied $r$ times approaches $\infty$ as more balls are thrown,
hence the time lag between two consecutive $r$-records is relatively small as compared with the total elapsed time. 
This covers the case of regularly varying $p_j$'s but excludes probabilities with exponential decay.

We  introduce a local time scale  to achieve that $r$-records occur at about constant rate, and obtain 
the Poisson approximation for the point processes of  $r$-record times. 
The joint convergence of records of different types  to independent Poisson processes is shown under the condition of regular  variation.

Similar Poisson approximation has been obtained
for the  familiar coupon-collectors problem~\cite{Ilienko}, which is the  occupancy scheme  with finitely many equiprobable  boxes.  Closely related work 
connecting record processes to the queueing theory
appeared  on the Ewens sampling model~\cite{GS}, where the probabilities $(p_j)$ are themselves random, chosen from the Poisson-Dirichlet/GEM distribution.
A characteristic property of the Ewens sampling model is that the indicators of records are independent \cite{Nacu}.

Following  Karlin's approach~\cite{Karlin} and the setting in much of the subsequent work, 
we first focus on a continuous time occupancy scheme where balls are thrown at epochs of a unit Poisson process.
This has the advantage over the discrete time scheme that the 
arrivals to distinct boxes occur  according to  {independent} Poisson processes. 
In Section \ref{depoisson}, we proceed with de-Poissonisation 
to obtain the Poisson approximation for records in the traditional model, where the discrete time variable coincides with the number of  balls thrown.
In Section \ref{discrete}, we build upon the exchangeability features   to develop a different  approach to the discrete time model. In the last section we consider occupancy  with random $(p_j)$, where  a mixed Poisson approximation to records is appropriate.

\section{Poissonised setup}

Suppose in the first instance that the balls labelled $1, 2,\ldots$ are thrown at epochs of a unit rate Poisson process $P=(P(t),~t\geq0)$.
 By the marking theorem,
box $j$ receives balls according to a Poisson process 
$P_j:=(P_j(t),~t\geq0), ~j\geq 1,$
with rate $p_j$, so the processes $P_j$ are independent and $P=\sum_{j\geq 1} P_j$.

Let $K_r(t):=|\{j: P_j(t)=r\}|, ~r\geq 1,$ denote the number of $r$-tons,  that is boxes  containing exactly $r$ balls at time $t$, and let 
$K(t):=\sum_{r\geq 1} K_r(t)$ be the number of boxes occupied by at least one ball.
The vector $(K_1(t), K_2(t),\ldots)$ encodes a random integer partition (possibly empty) induced by the allocation of the first  $P(t)$ balls.

We call {\it record} (time) any    jump time of the  process $K:=(K(t), ~t\geq 0)$.
For this and other nondecreasing counting processes, we shall use the common
convention to denote by the same symbol both the process  and the counting measure on Borel subsets of ${\mathbb R}_+$,
thus writing  
$K((u,t])=K(t)-K(u), 0\leq u<t$.
The counting process of $r$-tons $K_r:=(K_r(t),~t\geq0)$ has jumps $\pm 1$.
We  call {\it $r$-record}  (time)  any jump time  when 
$K_r$ increments by $+1$, hence  $K_{r-1}$ falls by $-1$ if $r>1$. Let 
$$
\beta_{jr}:=\inf\{t\geq 0: P_j(t)=r\},\quad j\geq 1,\quad r\geq 1.
$$
be the $r$-record time when box $j$ receives its $r$-th ball. 
These are well defined random variables, since 
the number of balls in each box grows to infinity, governed by the law of 
large numbers. For $j$ fixed,  the counting measure $\sum_{r\geq 1} \delta_{\beta_{jr}}$ is the Poisson random measure $P_j$ of arrivals to box $j$, where $\delta_t$ denotes the unit mass at $t$.

\begin{prop}\label{prop1}
For $\gamma\in[0,1]$ and $t\geq 0$, it holds that
\begin{eqnarray*}
{\mathbb E}[K(\gamma t)\,|\,K(t)]&\geq& \gamma K(t),\\
{\mathbb E}[K_r(\gamma t)\,|\,K_s(t)]&\geq& \binom{s}{r} \gamma^r(1-\gamma)^{s-r}K_s(t), ~~~s\geq r\geq 1.
\end{eqnarray*}
Therefore $(K(t)/t,~t\geq0)$ and
$(K_r(t)/t^r,~t\geq0)$ are reverse submartingales.
\end{prop}
\proof
 By the order statistic property of the Poisson process, given that at  time $t$ box $j$ has $s$  balls, the number of  arrivals 
to the box by time $\gamma t$ has the Binomial$(s,\gamma)$ distribution,
regardless of $p_j$. This implies the second inequality by the following estimates
\begin{align*}
{\mathbb E}[K_r(\gamma t)\,|\,K_s(t)]&=\sum_{j\geq 1}\mathbb{P}[P_j(\gamma t)=r\,|\,K_s(t)]\geq \sum_{j\geq 1}\mathbb{P}[P_j(\gamma t)=r,P_j(t)=s\,|\,K_s(t)]\\
&=\sum_{j\geq 1}\mathbb{P}[P_j(\gamma t)=r\,|\,P_j(t)=s]\mathbb{P}[P_j(t)=s\,|\,K_s(t)]\\
&=\binom{s}{r}\gamma^r(1-\gamma)^{s-r}\sum_{j\geq 1}\mathbb{P}[P_j(t)=s\,|\,K_s(t)]=\binom{s}{r}\gamma^r(1-\gamma)^{s-r}K_s(t).
\end{align*}
The first inequality follows along the same lines by noting that if a box contains $s$ balls at time $t$, then it was nonempty at time $\gamma t$ with probability at least $\gamma$, for every $s\geq 1$.
\endpf

We are interested in features of the point process of $r$-records
$B_r:=\sum_{j\geq 1} \delta_{\beta_{jr}}$.
Note that 
the records are the same as $1$-records, i.e. $B_1=K$, because  $K$ jumps when 
an empty box receives its first ball and becomes a singleton.

The sum of all $B_r$'s is the Poisson process $P$, but the processes $B_r$ themselves are not Poisson, rather possess a repulsion property known as the negative association; see~\cite{Last} for background. 

\begin{prop}\label{Pro1} Each record process $B_r, r\geq 1,$ is negatively associated.
\end{prop}

\proof The one-point process $\delta_{\beta_{jr}}$ is known to have the property of negative association. These are independent in $j$ by the independence of arrivals to boxes.
Hence $B_r$ is negatively associated as being a sum of the negatively associated processes.
\endpf
\noindent
In particular, the increments of $B_r$ over two disjoint intervals are negatively correlated.

The processes $B_r$ {and $B_s$, $r\neq s$,} are not independent. For instance, the first ball thrown after time $t=1$ is a $2$-record 
with zero probability if $\{B_1(1)=B_2(1)=0\}$,
and with positive probability if  $\{B_1(1)>0, B_2(2)=0\}$.

\section{The Bernstein function}

Without loss of generality, we assume that the boxes are labelled by decreasing popularity, that is $p_1\geq p_2\geq\cdots>0, ~\sum_{j\geq 1} p_j=1$. 
The probabilities $(p_j)$ are conveniently encoded into the infinite counting measure $\nu:=\sum_{j\geq1}\delta_{p_j}$ on $[0,1]$. 
This allows one to write sums over the boxes in the form of  integrals, 
$$\sum_{j\geq 1} f(p_j)=\int_{[0,\,1]} f(x)\nu({\rm d}x).$$

A dual way  to parametrise the model is to use
the {\it Bernstein function}
\begin{eqnarray}\label{BF}
\Phi(t):={\int_{[0,\,1]}} (1-e^{-tx})\nu({\rm d}x),\quad t\geq 0,
\end{eqnarray}
which uniquely determines $\nu$ and is important for the analysis. 
In the context of this integral representation, $\nu$ is sometimes called the  L\'{e}vy measure~\cite{Schilling}. 
The tilted  measure $\nu_1({\rm }x):=x\,\nu({\rm d}x)$ is normalised, and can be interpreted as the probability distribution of the popularity of the box hit by ball $1$.

The Bernstein function has a transparent 
 probabilistic meaning as the expected number of    boxes occupied by time $t$, that is 
$\Phi(t)={\mathbb E}[K(t)]$.  
Furthermore, the expected counts of $r$-tons
$\Phi_r(t):={\mathbb E}[K_r(t)]$ are expressible via the derivatives of $\Phi$ as
\begin{eqnarray*}
\Phi_r(t)&=&\frac{t^r (-1)^{r+1}}{r!}\Phi^{(r)}(t)=  \frac{t^r}{r!}{\int_{[0,1]}} e^{-tx}x^r \nu({\rm d}x).
\end{eqnarray*}
The formulas imply the recursion
\begin{equation}\label{der}
\Phi_{r+1}(t)=    \frac{r}{r+1}\Phi_r(t)-\frac{t}{r+1}\Phi_{r}'(t).
\end{equation}
Formulas for the variance and large-$t$ asymptotics are found in~\cite{GBarbour, Variance,  GHP,  Hwang, Karlin}.

In accord with  Proposition \ref{prop1} we have the monotonicity  $\Phi(t)/t\downarrow,\, \Phi_r(t)/t^r\downarrow$.
\begin{prop} The functions $\Phi_r$ satisfy for $0\leq\gamma\leq 1$ and $t\geq 0$ the  inequalities
\begin{eqnarray}
\label{est}
\Phi_r(\gamma t)&\geq & \binom{s}{r} \gamma^r (1-\gamma)^{s-r}\Phi_s(t),  ~1\leq r\leq s,\\
\label{est1}
\Phi_r(\gamma t)&\leq&   \gamma^r \left(\frac{ t^r  } {r!}  \right)^{1-\gamma}     \,\Phi_r^\gamma (t).
\end{eqnarray}
\end{prop}
\proof The first inequality follows from Proposition \ref{prop1}.
The second follows from Lemma \ref{L-app} in Appendix.
\endpf

\noindent
We note in passing that  a sharp constant in (\ref{est}) is obtained by  replacing the binomial coefficient  with the factor $\frac{s!}{r!} \left( \frac{e}{s-r}\right)^{s-r}$; see Eq.~(4.4) in~\cite{GBarbour}.

\begin{prop}\label{prop3} The point process of $r$-records has the intensity measure
$${\mathbb P}[B_r({\rm d}t)=1]=\frac{r\,\Phi_r(t)}{t} {\rm d}t,\quad t\geq0.$$
\end{prop}
\proof  By the order statistic property of the Poisson process, 
if box $j$ contains $r$ balls at time $t+{\rm d}t$, the  latest of them arrived within the interval $[t, t+{\rm d}t]$ with probability 
$r{\rm d}t/t$, regardless of $p_j$.
Given $K_r(t+{\rm d}t)=k$ the probability that a $r$-record occurs in the interval is $rk{\rm d}t/t$.
 The intensity formula follows by 
taking the expectation.
\endpf

The Bernstein function is concave, subadditive, satisfies $\Phi(t)\uparrow\infty$ and $\Phi(t)\ll t$ as $t\to\infty$, where $f(t)\ll g(t)$ means that $\lim_{t\to\infty}f(t)/g(t)=0$; see~\cite{Schilling}. In contrast,  the expected number of $r$-tons can be less regular.
Thus, as $t\to\infty$ the function $\Phi_r(t)$ may stay bounded, 
converge to $\infty$ or oscillate between a finite level and $\infty$; see~\cite{ GBarbour, Variance} 
for classification of the  modes of behaviour.
By  considering the $r$-records our first and foremost assumption will be  that $\Phi_r(t)\to\infty$, which is  equivalent to ${\rm Var}[K_r(t)]\to\infty$.
It is known, see~\cite{GBarbour}, that if $\Phi_r(t)\to\infty$ holds for some $r$ then also for all $r'\leq r$,  and a sufficient condition for  this is 
$p_j^r\ll \sum_{i\geq j+1} p_i^r, ~j\to\infty.$
This excludes the light-tailed distributions  $(p_j)$ like Poisson or geometric.

In the case $\limsup\Phi_r(t)<\infty$ an approximation to the process of $r$-records  can be sought on the contraction scale $\theta t, ~\theta\in[0,1]$,
similarly to~\cite{GS} or to the well known Poisson limit for record times in the extreme-value theory~\cite{Resnick}.
But this case falls outside the scope of the present note.


\section{Poisson approximation to \texorpdfstring{$r$}{r}-records}

Assuming $\Phi_r(t)\to\infty$, 
we aim at a local  Poisson approximation for  $B_r$.
The strategy is to fix some initial time $t_0$, called in the sequel {\it lower cutoff} and treated as a large parameter,  while introducing an auxiliary temporal variable $\theta>0$ to measure the size of a properly scaled time window.
These are related via  a time change
\begin{equation}\label{t-change}
t= t_0+ f(t_0)\,\theta,
\end{equation}
where $f$ is a suitable scaling function satisfying $f(t_0)\ll t_0$. Denote for shorthand $h:=f(t_0)\,\theta$.
The focus is then on
\begin{equation}\label{Q}
\widehat{B}_r(\theta):= B_r((t_0, t_0+h])=B_r((t_0, t_0+f(t_0)\,\theta]), 
\end{equation}
which is  the number of $r$-records arriving within the time window $h\ll t_0$.  
The point process $\widehat{B}_r$ as a component of $P$ is 
an instance of  `normalised ordered thinning'~\cite{Serfozo}.

By independence across the boxes, 
$\widehat{B}_r(\theta)$ has the Poisson-binomial distribution with   success probabilities 
\begin{equation}\label{q}
q_j:= \mathbb{P}[\beta_{jr}\in (t_0,t_0+h]]=\int_{t_0}^{t_0+h} e^{-tp_j} \frac{(tp_j)^{r-1}}{(r-1)!}\,p_j{\rm d}t,
\end{equation}
and expectation
\begin{equation}\label{meanQ}
\lambda_r:={\mathbb E}[\widehat{B}_r(\theta)]=\sum_{j\geq 1} q_j =\int_{t_0}^{t_0+h} \frac{r\Phi_r(t)}{t}{\rm d}t,
\end{equation}
as follows from Proposition \ref{prop3}. 

We proceed with the scaling function
\begin{equation}\label{marg-scale}
f(t_0)=\frac{t_0}{r\Phi_r(t_0)}.
\end{equation}
The value (\ref{marg-scale}) is the mean inter-arrival time in a homogeneous Poisson process with rate equal to the instantaneous rate of $B_r$ at time $t_0$.

Application of  Theorems 1.C(i) and 2.M from~\cite{Barbour}    yields the following estimate of the total variation distance to the Poiss$(\theta)$ distribution,
\begin{equation}\label{boundsTV}
d_{\rm TV}(\widehat{B}_r(\theta),{\rm Poiss}(\theta))\leq \frac{1-e^{-\lambda_r}}{\lambda_r}\,\sum_{j\geq 1} q_j^2+ |\lambda_r-\theta|. ~~~
\end{equation}
The first part is the seminal Barbour-Eagleson   bound on the total variation distance between 
$\widehat{B}_r(\theta)$ and ${\rm Poiss}(\lambda_r)$. The second part  
$|\lambda_r-\theta|$ is the bound for the distance between two Poisson distributions, see Lemma 1 in~\cite{Ruzankin}, which appears as the interpolation error caused
by  adopting parameter $\theta$ in place of the genuine mean $\lambda_r$.

To estimate the first part of the approximation error  (\ref{boundsTV})  start with the inequality
$$
q_j < \frac{e^{-t_0p_j} p_j^r}{(r-1)!}   \int_{t_0}^{t_0+h} t^{r-1}{\rm d}t< \frac{e^{-t_0p_j} p_j^r}{r!}\,  t_0^r \,\frac{ 2r h}{t_0}= \frac{e^{-t_0p_j} (p_j t_0)^r}{r!}  \,\frac{2 \theta}{\Phi_r(t_0)},
$$
which holds for all sufficiently large $t_0$.
Squaring this and summing over the boxes,  after some bookkeeping and application of 
(\ref{est}) to $\Phi_{2r}(2t_0)/\Phi_r(t)$ we obtain
$$
\sum_{j\geq 1} q_j^2<4 \theta^2  2^{-2r}{\binom{2r}{r}} \frac{\Phi_{2r}(2t_0)}{\Phi^2_r(t_0)}< \frac{4\theta^2}{\Phi_r(t_0)}.
$$
Quite satisfactorily, we see that $d_{\rm TV}(\widehat{B}_r(\theta), {\rm Poiss}(\lambda_r))$ tends to zero under the sole assumption 
$\Phi_r(t)\to\infty$ as $t\to\infty$. 

Exploring $|\lambda_r-\theta|$ turns more intricate.
From (\ref{meanQ}) and monotonicity
\begin{multline*}
\lambda_r\leq  \frac{r\Phi_r(t_0)}{t_0^r}\int_{t_0}^{t_0+h} t^{r-1}{\rm d} r = \frac{\Phi_r(t_0)}{t_0^r}\left( (t_0+h)^r-t_0^r\right)\\
=\Phi_r(t_0)\left(\frac{hr}{t_0} + \frac{r(r-1)h^2}{2t_0^2}+\cdots \right)\leq \theta + \frac{(r-1)}{2r}  \,\frac{\theta^2}{\Phi_r(t_0)}+\cdots,
\end{multline*}
where the remaining terms are of the smaller order as $t_0\to\infty$.
In the case $r=1$ the estimate simplifies as $\lambda_1\leq\theta$.
Similarly, a lower bound becomes
$$
\lambda_r\geq \frac{\Phi_r(t_0+h)}{(t_0+h)^r}\left( (t_0+h)^r     -t_0^r\right)=\frac{\Phi_r(t_0+h)}{\Phi_r(t_0)} \frac{t_0}{t_0+h}\,\theta - \frac{(r-1)}{2r}\frac{t_0^2}{(t_0+h)^2}\frac{\Phi_r(t_0+h)}{\Phi_r^2(t_0)}\,\theta^2+\cdots.
$$
It is seen that 
to achieve $|\lambda_r-\theta|\to 0$ it will be sufficient to fulfil
\begin{equation}\label{st-co}
\frac{\Phi_r(t_0+h)}{\Phi_r(t_0)}\to 1\quad {\rm as}\quad t_0\to\infty,\quad {\rm for}\quad h=\frac{t_0}{r\Phi_r(t_0)}\,\theta.
\end{equation}

Expanding
$\Phi_r(t_0+h)=\Phi_r(t_0)+\Phi_r'(t_0+uh)h, ~u\in[0,1],$
 after some plain manipulations using $\Phi_r(t)\to\infty$ and monotonicity, 
we show that a sufficient condition for (\ref{st-co})   to hold for each $\theta>0$   is
\begin{eqnarray}\label{wish2}
\frac{\Phi_{r+1}(t)}{\Phi_r^2(t)}\to 0\quad {\rm as}\quad t\to\infty.
\end{eqnarray}
In Section~\ref{sec:regularity}, we shall discuss regularity conditions that imply (\ref{wish2}).

The next result shows  that with  probability close to one there are no boxes that receive both $r$th and $(r+1)$st balls within the time window $[t_0,  t_0+h]$.
\begin{lemma}\label{no-doubles}
If $\Phi_r(t)\to\infty$ and either {\rm (\ref{st-co})} or {\rm (\ref{wish2})} holds, then 
$$
\sum_{j\geq 1} {\mathbb P}[t_0\leq \beta_{jr}<\beta_{j,r+1}\leq t_0+h]\to 0\quad {\rm as}\quad t_0\to\infty.
$$
\end{lemma}
\proof The generic term is the probability of the  event that box $j$ receives $r$th record and some other balls, hence 

\begin{multline*}
\sum_{j\geq 1} {\mathbb P}[t_0\leq \beta_{jr}<\beta_{j,r+1}\leq t_0+h]=\sum_{j\geq 1} \int_{t_0}^{t_0+h} e^{-tp_j}\frac{(tp_j)^{r-1}}{(r-1)!} p_j\left(1-e^{-(t_0+h-t)p_j}\right){\rm d}t\\
=\lambda_r-\sum_{j\geq 1}  \int_{t_0}^{t_0+h} e^{-(t_0+h)p_j}\frac{(tp_j)^{r-1}}{(r-1)!}p_j {\rm d}t=\lambda_r-  \frac{\Phi_r(t_0+h)}{ (t_0+h)^r} \left( (t_0+h)^r-{t_0}^r\right).
\end{multline*}
Here, both terms  converge to $\theta$ as before.
\endpf

We shall tacitly use two equivalent approaches to the functional convergence  of point processes on ${\mathbb R}_+:=[0,\infty)$; see Section 3 in~\cite{Resnick} or Section 11.1 in~\cite{Daley}. Recall that we identify $\widehat{B}_r$ with a random element of the space $M_p({\mathbb R}_+)$ of locally finite point measures on ${\mathbb R}_+$ by writing $\widehat{B}_r (A)=B_r(t_0+f(t_0)(A))$ for a Borel $A\subset{\mathbb R}_+$. We endow the space $M_p({\mathbb R}_+)$ with the topology of vague convergence. According to Lemma 11.1.XI in~\cite{Daley} the weak convergence on $M_p({\mathbb R}_+)$ is equivalent to the weak convergence of the corresponding cumulative processes in the Skorokhod space $D({\mathbb R}_+)$ of c\`{a}dl\`{a}g functions endowed with the $J_1$-topology. Thus, the random measures $\widehat{B}_r$ converge in distribution on $M_p({\mathbb R}_+)$ if and only if the random processes $(\widehat{B}_r(\theta), \theta\geq 0)$ converge in distribution on $D({\mathbb R}_+)$. Furthermore, by Theorem 11.1.VII in~\cite{Daley} both types of convergence are equivalent to the convergence of finite-dimensional distributions of $(\widehat{B}_r(\theta), \theta\geq 0)$ at the continuity points of the limit.

\begin{theorem}\label{T1}
If $\Phi_r(t)\to\infty$ as $t\to\infty$ and either of the  conditions {\rm (\ref{st-co})} or {\rm (\ref{wish2})} holds, then as $t_0\to\infty$ the process
$(\widehat{B}_r(\theta),~\theta\geq 0)$ converges in distribution to a Poisson process of unit rate. 
\end{theorem}
\proof {We first prove that for given $\theta>0$  the restriction of $\widehat{B}_r$ on $[0,\theta]$ converges to a Poisson process of unit rate on this interval. To this end, it suffices to show that the number and positions of atoms of $\widehat{B}_r$ on $[0,\theta]$ converge to the number and positions of atoms of the unit rate Poisson process on $[0,\theta]$. The number of points of $\widehat{B}_r$ on $[0,\theta]$ converges to ${\rm Poiss}(\theta)$, because both parts of the error (\ref{boundsTV}) approach $0$. For the convergence of positions recall a} familiar re-statement of the order statistic property:  the Poisson process on $[0,\theta]$ can be characterised as a mixed binomial process, whose number of points has  ${\rm Poiss}(\theta)$ distribution; see Section 3 in~\cite{Kallenberg}. Let $A$ be the event that no box receives the $r$th and $(r+1)$st balls within the same time window $[t_0, t_0+h]$.
Given $A$, the process $\widehat{B}_r$ restricted to ${[0,\theta]}$ is mixed binomial, because a single arrival  to a box within the time window is uniformly distributed there
and under our scaling (\ref{t-change}), (\ref{marg-scale}) the arrival becomes uniformly distributed on $[0,\theta]$.
By Lemma \ref{no-doubles} ${\mathbb P}[A]\to 1$ as $t_0\to\infty$, which together with the convergence of $\widehat{B}_r(\theta)$ implies convergence to the Poisson process on $[0,\theta]$. Since this holds for every $\theta$, the proof is complete.
\endpf

\section{A conditional limit theorem}

So far we have been concerned with the features of the process of $r$-records averaged over the history prior to the cutoff time $t_0$.
However, the online observer being aware of the history will see different patterns depending on the allocation at time $t_0$. 
Two extreme cases are accumulation of all balls in just one box and the situation when each arrival before $t_0$ hits a different box.
In this section we show that typically the past allocation does not impact the future, that is with high probability the same Poisson approximation works even conditionally on the full history
of the occupancy process.

{Let $\mathcal{F}_{t_0}:=\sigma(P_j(t),j\geq 1,t\leq t_0)$ be the sigma-algebra generated by the occupancy counts up to and including time $t_0$. 
\begin{prop}
Under the assumptions of {\rm Theorem~\ref{T1}}, as $t_0\to\infty$  the conditional distribution of $(\widehat{B}_r(\theta),\theta\geq 0)$ given $\mathcal{F}_{t_0}$ converges in probability to the distribution of a Poisson process of unit rate.
\end{prop}
\proof
According to Lemma~\ref{no-doubles}
$$
\widehat{B}_r(\theta)=\sum_{j\geq 1}\mathbbm{1}_{\{\beta_{jr}\in (t_0,t_0+f(t_0)\theta\}}=\sum_{j\geq 1}\mathbbm{1}_{\{P_j(t_0+f(t_0)\theta)\geq r,P_j(t_0)<r\}}=\sum_{j\geq 1}\mathbbm{1}_{\{P_j(t_0+f(t_0)\theta)\geq r,P_j(t_0)=r-1\}}+\varepsilon,
$$
where $\varepsilon=\varepsilon(\theta)$ converges to zero in probability locally uniformly with respect to $\theta$, as $t_0\to\infty$. Denote by $P_j^{\prime}(t):=P_j(t+t_0)-P_j(t_0)$, $j\geq 1$, $t\geq 0$ and note that $(P_j^{\prime},j\geq 1)$ is independent of $\mathcal{F}_{t_0}$. Thus we need to check that the conditional distribution of the process
\begin{equation}\label{fixedTheta}
\theta\mapsto \sum_{j:P_j(t_0)=r-1}\mathbbm{1}_{\{P^{\prime}_j(f(t_0)\theta)\geq 1\}},
\end{equation}
given $\mathcal{F}_{t_0}$, converges to the desired Poisson limit. Note that for fixed $\theta$, given $\mathcal{F}_{t_0}$, 
the variable (\ref{fixedTheta})
is a sum of independent Bernoulli variables. 
Arguing in the same way as in the proof of Theorem~\ref{T1}, we see that it suffices to prove that the conditional distribution of (\ref{fixedTheta})
converges in probability to ${\rm Poiss}(\theta)$, for every fixed $\theta>0$. Put
$$
E_{t_0,r}:=\mathbb{E}\left[\sum_{j:P_j(t_0)=r-1}\mathbbm{1}_{\{P^{\prime}_j(f(t_0)\theta)\geq 1\}}\Big| \mathcal{F}_{t_0}\right]=\sum_{j\geq 1}(1-e^{-p_j f(t_0)\theta})\mathbbm{1}_{\{P_j(t_0)=r-1\}},
$$
and note that
$$
\mathbb{E} [E_{t_0,r}]=\sum_{j\geq 1}(1-e^{-p_j f(t_0)\theta})e^{-p_j t_0} \frac{(p_j t_0)^{r-1}}{(r-1)!}=\theta-\sum_{j\geq 1}(p_jf(t_0)\theta-1+e^{-p_j f(t_0)\theta})e^{-p_j t_0} \frac{(p_j t_0)^{r-1}}{(r-1)!}.
$$
Using an elementary inequality $x-1+e^{-x}\leq x(1-e^{-x})$, $x\geq 0$, we see that
\begin{align}\label{eq:conditional_proof1}
0&\leq\sum_{j\geq 1}(p_jf(t_0)\theta-1+e^{-p_j f(t_0)\theta})e^{-p_j t_0} \frac{(p_j t_0)^{r-1}}{(r-1)!}\notag\\
&\leq\sum_{j\geq 1}p_jf(t_0)\theta(1-e^{-p_j f(t_0)\theta})e^{-p_j t_0} \frac{(p_j t_0)^{r-1}}{(r-1)!}\notag\\
&= \frac{rf(t_0)\theta}{t_0}\left(\Phi_r(t_0)-\frac{t_0^r}{(t_0+f(t_0)\theta)^r}\Phi_r(t_0+f(t_0)\theta)\right)\\
&=1-\frac{t_0^r}{(t_0+f(t_0)\theta)^r}\frac{\Phi_r(t_0+f(t_0)\theta)}{\Phi_r(t_0)}\to 0\notag,
\end{align}
by the assumption~\eqref{st-co}. Similarly, using $(1-e^{-x})^2\leq x(1-e^{-x})$, $x\geq 0$, we obtain from~\eqref{eq:conditional_proof1}
\begin{multline}\label{eq:conditional_proof2}
{\rm Var}[E_{t_0,r}] \leq \sum_{j\geq 1}(1-e^{-p_j f(t_0)\theta})^2 e^{-p_j t_0} \frac{(p_j t_0)^{r-1}}{(r-1)!}\\
\leq \sum_{j\geq 1}p_jf(t_0)\theta(1-e^{-p_j f(t_0)\theta})e^{-p_j t_0} \frac{(p_j t_0)^{r-1}}{(r-1)!}\to 0.
\end{multline}
Thus, $E_{t_0,r}$ converges to $\theta$ in probability as $t_0\to\infty$ by Chebyshev's inequality. In view of~\eqref{boundsTV}, given $\mathcal{F}_{t_0}$,
$$
d_{\rm TV}\left(\sum_{j:P_j(t_0)=r-1}\mathbbm{1}_{\{P^{\prime}_j(f(t_0)\theta)\geq 1\}},{\rm Poiss}(\theta)\right)\leq\frac{1-e^{-E_{t_0,r}}}{E_{t_0,r}}\sum_{j\geq 1}(1-e^{-p_j f(t_0)\theta})^2\mathbbm{1}_{\{P_j(t_0)=r-1\}}+|E_{t_0,r}-\theta|.
$$
The right-hand side converges to zero in probability by Markov's inequality, since
$$
\mathbb{E}\left[\sum_{j\geq 1}(1-e^{-p_j f(t_0)\theta})^2\mathbbm{1}_{\{P_j(t_0)=r-1\}}\right]=\sum_{j\geq 1}(1-e^{-p_j f(t_0)\theta})^2 e^{-p_j t_0} \frac{(p_j t_0)^{r-1}}{(r-1)!}\to 0,
$$
as in~\eqref{eq:conditional_proof2}.
\endpf

\section{Regularity and growth}\label{sec:regularity}

If $\Phi_r'(t)\geq 0$ then $\Phi_{r+1}(t)\leq \frac{r}{r+1}\Phi_r(t)$ by (\ref{der}), hence (\ref{wish2}) is always true
if $t\to \infty$ along a sequence of increase points of $\Phi_r$ (subject to the only condition $\Phi_r(t)\to\infty$).
However,  we could neither verify (\ref{wish2}) in full generality  nor construct a counter-example in terms of $(p_j)$ or $\Phi$.
In this section  we give various conditions to
ensure (\ref{wish2}) or directly (\ref{st-co}), hence the Poisson convergence of  the process of $r$-records, in accord with  Theorem \ref{T1}.

\subsection{Regular variation}\label{secRV}
Karlin's condition of regular variation~\cite{Karlin} reads as
\begin{equation}\label{RV}
\nu[x,1]\sim x^{-\alpha}\ell(1/x),\quad x\to 0+,
\end{equation}
where $0\leq \alpha\leq 1$  is  the {\it index} and $\ell$ is some  function slowly  varying at $\infty$.
In the {\it proper} case $0<\alpha<1$ this implies that,  asymptotically,  the expectations of the $r$-ton counts  only differ by constant factors:
\begin{equation}\label{Phi-r}
 \Phi(t)\sim \Gamma(1-\alpha)\, t^\alpha\ell(t)\quad {\rm and}\quad \Phi_r(t)\sim \frac{\alpha\Gamma(r-\alpha)}{r!\Gamma(1-\alpha)}\,\Phi(t).
\end{equation}

In the case of {\it rapid} variation $\alpha=1$ the slowly varying factor must satisfy $\ell(t)\to 0$ as $t\to\infty$ (to agree with $\sum_{j\geq 1} p_j=1$), in which case the asymptotic formulas are
\begin{equation}\label{rapid-Phi}
\Phi(t)\sim t\ell_1(t),\quad \Phi_1(t)\sim\Phi(t)\quad {\rm and~}\quad \Phi_r(t)\sim \frac{1}{r(r-1)}t\ell(t)\quad {\rm for}\quad r\geq 2,
\end{equation}
where $\ell_1\gg \ell$ is another slowly varying function, thus $\Phi_1\gg \Phi_r$ for $r\geq 2$.

Speaking of the case of {\it slow} variation $\alpha=0$ we shall mean a slightly  stronger condition 
\begin{equation}\label{slow}
\nu_1[0,x]\sim x\ell_0(1/x),\quad x\to 0+,
\end{equation}
with slowly varying $\ell_0(t)\to\infty$, $t\to\infty$; then (\ref{RV}) holds with some $\ell\gg\ell_0$, and 
\begin{equation}\label{slow-Phi}
\Phi_r(t)\sim\frac{1}{r}\ell_0(t),\quad \Phi(t)\sim \ell(t),\quad t\to\infty,
\end{equation}
so the $\Phi_r$'s are of the same order while $\Phi(t)\gg \Phi_r(t)$ for all  $r\geq 1$. Note that the mean number $r\Phi_r(t)$ of balls contained in $r$-ton boxes is asymptotic to the same function $\ell_0(t)$
regardless of $r$.
See~\cite{GHP} for conditions of slow variation expressed in terms of $(p_j)$ or $\nu$.

\subsection{A weaker form of regular variation}
Under the regular variation (understood as (\ref{slow}) if $\alpha=0$) we obviously have
\begin{equation}\label{wish1}
\liminf\limits_{\gamma\downarrow 1}
\liminf\limits_{t\to\infty} \frac{\Phi_r(\gamma t)}{\Phi_r(t)}\geq 1,
\end{equation}
which  in turn implies (since $\Phi_r(t)/t^r$  decreases) the desired (\ref{st-co}) for all $\theta>0$. Condition (\ref{wish1}) itself is well known in the Tauberian theory; see p.~19 in~\cite{Goldie}. If it holds  we have 
$$
\lim\limits_{\gamma\to1}
\limsup\limits_{t\to\infty} |{\Phi_r(\gamma t)}/{\Phi_r(t)}-1|=0.
$$
The latter asymptotic condition defines the class of  {\it pseudo-regularly varying}
functions, treated in a  recent monograph~\cite{Klesov}.

\subsection{Slow decrease}\label{sl-d}

Assuming $\Phi_r(t)\to\infty$, it is easy to check that
\begin{eqnarray}\label{wish3}
\frac{\Phi_{r+1}(t)}{\Phi_r^2(t)}\to 0\quad \Longleftrightarrow\quad \left(\frac{t}{\Phi_r(t)}\right)'\to 0.
\end{eqnarray}
Writing
$$
\frac{1}{\Phi_r(t)} =\frac{1}{t} \int_a^t   \left(\frac{u}{\Phi_r(u)}\right)' {\rm d}u+o(1),
$$
we have the left-hand side converging to $0$, which only forces the Ces{\'a}ro summability 
hence a priori does not  guarantee  convergence of the integrand. 
Thus some extra assumption to limit the variability of $\Phi_r$ seems inevitable.

In the classic text \cite[Section 6.2]{Hardy} a function $g$ is called {\it slowly decreasing} if 
$$
\liminf (g(\gamma t)-g(t))\geq 0
$$ 
for $t\to\infty$, $\gamma\geq 1$ and $\gamma\to 1$. A sufficient condition for this is $g'(t)>-c/t$.
Applying~\cite{Hardy} Theorem 68, condition (\ref{wish2}) holds if any of the involved  functions $\Phi_r$ and $\Phi_{r+1}$  is slowly decreasing.

\subsection{A bounded ratio}

Given $\Phi_r(t)\to\infty$, a sufficient condition for  (\ref{wish2}) is
\begin{equation}\label{rat}
\limsup\limits_{t\to\infty} \frac{\Phi_{r+1}(t)}{\Phi_r(t)}<\infty.
\end{equation}
Now suppose that for some $1<\gamma_1<\gamma_2$
\begin{equation}\label{quot}
\liminf\limits_{t\to\infty} \sup\limits_{\gamma \in [\gamma_1,\gamma_2]}\frac{\Phi_{r+1}(\gamma t) }{\Phi_{r+1}( t) }>0.
\end{equation}
Then from (\ref{est}) for  any such $\gamma$
$$
\frac{\Phi_{r+1}(\gamma t)}{\Phi_r(t)}\leq \frac{\gamma^{r+1}}{(r+1)(\gamma-1)}\leq \frac{\gamma_2^{r+1}}{(r+1)(\gamma_1-1)},
$$
whence choosing $\varepsilon$ less than the limes inferior  in (\ref{quot}), for large enough $t$
$$
\frac{\Phi_{r+1}(t)}{\Phi_r(t)}< \frac{\gamma_2^{r+1}}{\varepsilon (r+1)(\gamma_1-1)},
$$
so (\ref{rat}) is implied by  (\ref{quot}). By the same token, (\ref{rat}) also follows if  condition (\ref{quot}) is  imposed on $\Phi_r$ instead of $\Phi_{r+1}$.

\subsection{Superlogarithmic growth}

Finally, we show that the superlogarithmic growth condition 
\begin{equation}\label{suplog}
\Phi_r(t)\gg \log t
\end{equation}
is sufficient for (\ref{st-co}). Indeed, from (\ref{est}) and (\ref{est1}) we may estimate (\ref{wish2}) as

\begin{align*}
\frac{\Phi_{r+1}(t)}{\Phi^2_r(t)}&\leq \frac{\Phi_r(\gamma t)}{(r+1)\gamma^r (1-\gamma)\Phi_r^2(t)}\\
&\leq \frac{2\, t^{r(1-\gamma)}\,\Phi_r^\gamma(t)}  {(r!)^{1-\gamma} (r+1)(1-\gamma)\,\Phi_r^2(t)}\leq \frac{ t^{r(1-\gamma)}}   {(1-\gamma) \Phi_r^{2-\gamma}(t)}\,.
\end{align*}
If $\Phi_r(t)= g(t)\log t$  with some $g(t)\to\infty$, with the 
choice of the   parameter  $\gamma=1-(\log t)^{-1}$ the estimate becomes
$$
\frac{\Phi_{r+1}(t)}{\Phi^2_r(t)}=
\frac{e^r}{g(t)}
 \{g(t)\log t\}^{-1/\log t}\to 0.
$$

A similar check shows that the bound does not converge to $0$ whenever $\limsup g(t)<\infty$.  This only means tha the bound is inconclusive, because
  $\Phi_{r+1}(t)/\Phi^2_r(t)\to 0$ may still hold for $\Phi_r$ arbitrarily slowly growing, e.g. under the condition  
of slow variation (\ref{slow}).

In the view of discussion in Section \ref{sl-d}, it looks unexpected  that the growth rate (\ref{suplog}) assumes the role of a Tauberian condition.

\section{Multivariate processes of records}

For the joint convergence of $r$-record processes to Poisson processes with constant rates, 
one needs a common scaling function $f$ in (\ref{t-change}) to serve different types of records. 
But  this is only possible under a condition of regular variation.  Indeed,
note that $L_r(t):=\Phi_r(t)/(t^r/r!)$  is the Laplace transform of the measure $x^r\nu({\rm d}x)$ on $[0,1]$, hence {relation~\eqref{der}} between $\Phi_r$ and $\Phi_{r+1}$ becomes
$$L_r(t)=\int_t^\infty L_{r+1}(u){\rm d}u,$$
and therefore
$$\frac{\Phi_{r+1}(t)}{\Phi_r(t)}=\dfrac{t L_{r+1}(t)}{(r+1)\int_t^\infty L_{r+1}(u){\rm d}u}.$$
If  the ratio converges as $t\to\infty$ to some constant, say $(r-\alpha)/(r+1)$, then by Karamata's theorem, see Theorem 1.6.1 in~\cite{Goldie}, $L_{r+1}$ is regularly varying with index $-(r+1)+\alpha$, thus both functions $\Phi_{r+1}$ and $\Phi_r$ are regularly varying with index $\alpha$.

\begin{theorem}\label{T2}
Suppose $\Phi$ is regularly varying with index $0\leq\alpha\leq1$.
Then with  the scaling function $f(t_0)$
the processes $\widehat{B}_r$ converge jointly  in distribution, as $t_0\to\infty$,  to independent homogeneous Poisson processes with rates 
$\rho_r$, where
$$
\begin{array}{llll}
{\rm for}~0<\alpha<1: &
 f(t_0)=\dfrac{t_0}{\Phi(t_0)}~~{\rm and~} &
\rho_r = \dfrac{\alpha \Gamma(r-\alpha)}{(r-1)!\, \Gamma(1-\alpha)}, & r\geq1;\\
{\rm for}~\alpha=1:&
 f(t_0)=\dfrac{t_0}{\Phi_2(t_0)}~~{\rm and~}&
\rho_r  =\dfrac{2}{r-1}, &r\geq2;\\
{\rm for}~\alpha=0~{\rm under~(\ref{slow})}:&
 f(t_0)=\dfrac{t_0}{\Phi_1(t_0)}~~{\rm and~} &
\rho_r =1, &r\geq1.
\end{array}
$$
\end{theorem}
\proof
The rates  are clear from (\ref{Phi-r}), (\ref{rapid-Phi}) and (\ref{slow-Phi}).
The marginal convergence of the scaled $r$-record processes follows from Theorem \ref{T1}.
For the joint convergence, we need  to show that for each $r$ the processes
$\widehat{B}_1,\ldots, \widehat{B}_r$ converge jointly to  a multivariate Poisson process.
To that end, we apply Corollary 11.2.VII in~\cite{Daley}. With the marginal convergence at hand, the convergence of intensity measures, see Eq.~(11.2.11) in~\cite{Daley}, holds automatically. It remains  to verify the condition
$$
\sum_j {\mathbb P}\left[\sum_{s=1}^r \delta_{\beta_{js}}([t_0, t_0+f(t_0)\theta])\geq 2\right]\to 0\quad {\rm as}\quad  t_0\to\infty,
$$
which is Eq.~(11.2.10) in~\cite{Daley}. But this is satisfied  by Lemma \ref{no-doubles}, since 
$$
\sum_{s=1}^r \delta_{\beta_{js}}([t_0, t_0+h])\geq 2
$$ 
means precisely that $t_0\leq \beta_{js}<\beta_{j,s+1}\leq t_0+h$, for some $1\leq s<r$.
\endpf
\noindent
A multivariate counterpart of  (\ref{boundsTV}) can be derived from estimates in~\cite{Roos}.

In the case $\alpha=1$ of rapid variation the scaling in Theorem \ref{T2}   is not suitable for the process $B_1$, because the  $1$-records are then much more frequent than records of any other type $r>1$.

\endpf

\section{The discrete-time model: de-Poissonisation}\label{depoisson}

We turn next to  the occupancy scheme, where balls are thrown at  discrete times.
De-Poissonisation is a folk name for   methods aiming to derive properties of `fixed-$n$' models from their `Poiss$(n)$' counterparts.  
For $K(t), K_r(t)$ the de-Poissonisation relies on concentration properties of the Poisson distribution, which enable efficient coupling where the values of the variables in both models are determined by much the same bulk of balls; see, for example, Section 6.2 in~\cite{Moderate}. For $r$-record processes such coupling can be constructed  if $\Phi_r(t)\ll t^{1/2}$, but is impossible in models with  $\Phi_r(t)\gg t^{1/2}$ where the time window has the order  smaller than the square root fluctuations. Fortunately, if $\Phi_r(t)\to\infty$ the de-Poissonisation can be universally justified due to the square root insensitivity  to the choice of  the lower cutoff $t_0$.

Denote $\mathcal{K}_{n}$ the total number of occupied boxes when $n$ balls are allocated and   $\mathcal{K}_{n,r}$ the number of boxes occupied by exactly $r$ balls, $r\geq 1$. Thus,
$$
K(t)=\mathcal{K}_{P(t)},\quad K_r(t)=\mathcal{K}_{P(t),r},\quad r\geq 1.
$$
For the moments we shall use the approximate formulas
\begin{equation}\label{disc-mean}
{\mathbb E}[\mathcal{K}_{n}]=\Phi(n)+o(1),~~~  {\mathbb E}[\mathcal{K}_{n,r}]=\Phi_r(n)+o(1),~~~~~n\to\infty,
\end{equation}
found in~\cite{Variance} along with an explicit  fixed-$n$ analogue of (\ref{BF}).

We shall also need a vector of occupancy counts $(\mathcal{P}_{j}(n),j\geq 1)$, where $\mathcal{P}_j({n})$ is the number of balls in  box $j$ after allocating $n$ balls. A discrete-time counterpart of the record time $\beta_{j r}$ is
$$
\widetilde{\beta}_{j r}:=\min\{n\in\mathbb{N}:\mathcal{P}_j(n)=r\},\quad j\geq 1,\quad r\geq 1,
$$
and we have an obvious relation $\widetilde{\beta}_{j r}=P(\beta_{j r})$. Therefore, the point process of discrete-time $r$-records 
$$
\mathcal{B}_r:=\sum_{j\geq 1}\delta_{\widetilde{\beta}_{j r}},\quad r\geq 1,
$$
is the random measure on $\mathbb{N}$, which may be represented as the push-forward of $B_r$ under the random mapping $ t\mapsto P(t)$, that is $\mathcal{B}_r=B_r\circ P^{-1}$. Note that $B_r$ and $P$ in this representation  are dependent making analysis of $\mathcal{B}_r$ harder. Let  $\widehat{\mathcal{B}}_r(\theta):=\mathcal{B}_r((n_0, n_0+f(n_0)\,\theta]), ~\theta\geq0$.

\begin{theorem}\label{T3} 
Fix $r\geq 1$.
If $\Phi_r(t)\to\infty$ as $t\to\infty$ and either of the conditions {\rm (\ref{st-co})} or {\rm (\ref{wish2})}  is satisfied, then with the scaling function 
$f(n_0)= {n_0}/({r\Phi_r(n_0)}),$
as $n_0\to\infty$  the process
$\widehat{\mathcal{B}}_r$ converges in distribution  to a Poisson process of unit rate.
\end{theorem}
\vskip-0.15cm
To prove the theorem we need two auxiliary lemmas. Recall that we assume $p_1\geq p_2\geq \cdots$. For $t>1$, let $j(t)\in\mathbb{N}$ be the unique index such that
$$
p_{j(t)}> \frac{2\log t}{t}\geq p_{j(t)+1}.
$$
We regard a box with $j\leq j(t_0)$ as  `popular' since
 for large time $t_0=n_0$, 
it is likely to contain more than $r$ balls, both in Poisson and discrete time occupancy schemes.
The intuition suggests that the popular boxes make negligible contribution to the normalised record processes.
We prefer to justify this in the {continuous} time setting, the other case being completely analogous.
\begin{lemma}\label{lemma-large} For $r\geq 1$, 
$$
{\mathbb P}[P_j(t)\leq r~{\rm for~some~}j\leq j(t)]\to 0\quad {\rm as}\quad t\to\infty.
$$
\end{lemma}
\proof For $j\leq j(t)$ we have  ${\mathbb E}[P_j(t)]\geq 2\log t$. The assertion follows from the elementary estimate 
\begin{equation}\label{eq:poisson_left_tail}
{\mathbb P}[P_j(t)\leq r]\leq {\frac{(r+1)}{r!}}\frac{ (2\log t)^r}{t^2}\quad
{\rm  for}\quad t>e^{r/2},
\end{equation}
by observing that $j(t)\leq t/(2\log t)$ because there are at most $1/p$ boxes with probability larger $p$. {Bound~\eqref{eq:poisson_left_tail} is a consequence of the chain of estimates: for $r\leq 2\log t\leq \theta$,
$$
{\mathbb P}[{\rm Poiss}(\theta)\leq r]\leq (r+1){\mathbb P}[{\rm Poiss}(\theta)=r]= \frac{(r+1)}{r!}e^{-\theta}\theta^r\leq \frac{(r+1)}{r!}e^{-2\log t}(2\log t)^r,
$$
where for the last inequality we used that $\theta\mapsto e^{-\theta}\theta^r$ is decreasing for $\theta\geq r$.}
\endpf

\noindent
By the lemma, a truncated version of the normalised $r$-record process, 
\begin{equation}\label{truncated}
\theta\mapsto \sum_{j\geq j(t_0)}\delta_{\beta_{j r}}((t_0,t_0+f(t_0)\,\theta]),\quad \theta\geq 0,
\end{equation}
for large $t_0$ coincides with $\widehat{B}_r$ with probability close to one.

For $n=1,2,\ldots$ the discrete time allocations are naturally identified with the configuration of balls in boxes at random times $S_n$, where
$S_n:=\min\{t\geq 0:~P(t)=n\}$ is the $n$-th arrival time in the Poisson process $P$.
Consider the process
$$
\widehat{B}_r^*(\theta):=
\sum_{j \geq 1}\delta_{\beta_{j r}}((S_{n_0},  S_{n_0}+f(n_0)\theta]),~~\theta\geq0,
$$
which has the same window size as $\widehat{B}_r(\theta)$ but the lower cutoff $t_0=n_0$ is replaced by the $n_0$th point of $P$.
Replacing $n_0$ by $S_{n_0}$ is a nontrivial step which turns possible in full generality due to a key
observation from~\cite{GBarbour} that the counts of balls within unpopular boxes at large times $t_0=n_0$ are similar for both Poisson and discrete-time schemes.

\begin{lemma}\label{lem:shift_of_window}
Under conditions of {\rm Theorem \ref{T1}}, 
for $r\geq 1$, as $n_0\to\infty$ the process
$(\widehat{B}_r^*(\theta), ~\theta\geq0)$
converges in distribution to a Poisson process with unit rate.
\end{lemma}
\proof
Lemma \ref{lemma-large} implies that it is sufficient to restrict summation over $j>j(n_0)$. 
Since $P(S_{n_0})=n_0$, we may apply the total variation
 estimate~(2.6) from~\cite{GBarbour}, which in our notation reads as
\begin{multline}\label{eq:GneBar_TV}
d_{\rm TV}\left((\mathcal{P}_j(n_0),~j> j(n_0),~(P_j(n_0),~j> j(n_0)\right)\\
=d_{\rm TV}\left((P_j(S_{n_0}),~j> j(n_0)), (P_j(n_0),~j> j(n_0)\right)\leq \sum_{j>j({n_0})} p_j.
\end{multline}
Letting $n_0\to\infty$ gives $j(n_0)\to\infty$ hence the right side approaches $0$.
That is to say, the occupancy numbers in unpopular boxes are likely to be the same at times $n_0$ and $S_{n_0}$.
On the other hand, 
 the $r$-record process after time $n_0$ depends on the history only through the allocation of balls at this time. The assertion now follows
from the convergence of $\widehat{B}_r$.
\endpf


We are now in position to prove Theorem~\ref{T3}.

\proof {The $r$-record process is nondecreasing, therefore the result will follow from Lemma \ref{lem:shift_of_window} by a sandwich argument provided we can justify that
$${\mathbb P}[ \widehat{B}_r^*(\theta_i-\varepsilon)\leq \widehat{B}_r(\theta_i)\leq \widehat{B}_r^*(\theta_i+\varepsilon),i=1,\ldots,m]\to 1,~~~n_0\to\infty,$$
for every fixed $m\in\mathbb{N}$, $0\leq \theta_1<\cdots< \theta_m$ and $\varepsilon>0$. Clearly, it suffices to consider the case $m=1$. Note that the event 
$$
\Big\{P\big((S_{n_0},S_{n_0}+f(n_0)(\theta_1-\varepsilon)]\big)\subseteq (n_0,n_0+f(n_0)\theta_1]\subseteq P\big((S_{n_0},S_{n_0}+f(n_0)(\theta_1+\varepsilon)]\big)\Big\}
$$
implies the event $\{\widehat{B}_r^*(\theta_1-\varepsilon)\leq \widehat{B}_r(\theta_1)\leq \widehat{B}_r^*(\theta_1+\varepsilon)\}$. But the former, by monotonicity of $P$ and the fact that $P(S_{n_0})=n_0$, coincides with the event
$$
\Big\{P(S_{n_0}+f(n_0)(\theta_1-\varepsilon))-P(S_{n_0})\leq f(n_0)\theta_1\leq P(S_{n_0}+f(n_0)(\theta_1+\varepsilon))-P(S_{n_0})\Big\}.
$$
The probability of this event tends to one by the weak law of large numbers for $(P(S_{n_0}+t)-P(S_{n_0}),t\geq 0)$ which is again a unit rate Poisson process by the strong Markov property of $P$.}
\endpf

A discrete-time version of Theorem~\ref{T2} also holds true. The proof proceeds along the same lines, does not involve any new ideas and is therefore omitted. 

\begin{theorem}\label{T4}
Under the assumptions of {\rm Theorem~\ref{T2}} the random  processes $\widehat{\mathcal{B}}_r$ converge jointly in distribution, as $n_0\to\infty$,  to independent homogeneous Poisson processes with rates $\rho_r$.
\end{theorem}

\section{The discrete-time model: the use of exchangeability}\label{discrete}

Records arriving at large times typically emerge due to unpopular  boxes  that rarely change their occupancy status. Therefore at later stages all relatively recent $r$-record balls, for fixed $r$,  are likely to belong to the $r$-ton boxes.  This intuitive feature combined with  the intrinsic exchangeability of the occupancy scheme  leads to another approach to the Poisson approximation, which we 
sketch in the discrete time setting and under the assumption of regular variation. The aim here is to show weak convergence to a multivariate Poisson distribution of the random vector or record counts $(\mathcal{B}_s((n_0, n_1]), ~1\leq s\leq r)$
for suitable time window $(n_0, n_1]$ with $n_0\to\infty$. For simplicity we also exclude the case $\alpha=1$ where the multivariate approximation holds for $r\geq 2$. Under multivariate Poisson distribution we understand the distribution of an integer  vector with independent univariate Poisson components.

For given $n_1$ define $\mathcal{C}_r$ to be the random set of balls contained in the $r$-ton boxes present at time $n_1$. Formally, 
$\mathcal{C}_r$ is
 a random point measure on $\{1,\ldots, n_1\}$, 
though we do not include $n_1$ in the notation.
By the definition, $n$ is an atom of
$\mathcal{C}_r$ if and only if $n=\widetilde{\beta}_{jr}\leq  n_1<\widetilde{\beta}_{j,r+1}$ for some $j$, and by the virtue of
 (\ref{disc-mean})  this event has probability
$$
\frac{r}{n_1}\,{\mathbb E}[\mathcal{K}_{n_1,r}]=\frac{r}{n_1}\,\Phi_r(n_1)+o\left(\frac{1}{n_1}\right),
$$
which is the same for all $n\leq n_1$.
To compare $\mathcal{C}_r$ with $r$-record counts
in terms of their means choose  $n_0<n_1$ and observe that
\begin{eqnarray}
\label{mean1}
{\mathbb E}[\mathcal{C}_r((n_0,n_1])]&=&\frac{r(n_1-n_0)}{n_1}\,\Phi_r(n_1)+ o\left(\frac{n_1-n_0}{n_1}\right),\\
\label{mean2}
{\mathbb E}[\mathcal{B}_r((n_0,n_1])]&=&\sum_{n=n_0+1}^{n_1} \frac{r\Phi_r(n)}{n}+o\left(\frac{n_1-n_0}{n_1}\right).
\end{eqnarray}
Setting $n_1=n_0+\lfloor f(n_0)\theta\rfloor$ with the scaling function as in  (\ref{marg-scale}),   we obtain under  conditions
of Theorem \ref{T1} that
$
{\mathbb E}[\mathcal{B}_r((n_0,n_1])- \mathcal{C}_r((n_0,n_1])]\to 0
$ 
as $n_0\to\infty$ locally uniformly in $\theta\geq 0$.

The multivariate Poisson approximation to  records will be justified in several steps.

\vspace{2mm}
\noindent
{\bf Step 1:} Approximation of $(\mathcal{B}_s((n_0, n_1]), ~1\leq s\leq r)$ by $(\mathcal{C}_s((n_0, n_1]), ~1\leq s\leq r)$.

\vspace{2mm}
The point process $\mathcal{C}_r$ is much better tractable than $\mathcal{B}_r$, because $\mathcal{C}_r$ is exchangeable, that is invariant under re-labelling of balls $1,\ldots,n_1$ by permutations. If a generic ball $n$ at time $n_1$ belongs to an $s$-ton for some $s\leq r$, then $n$ is also a $s_1$-record time for some $s_1\leq s$. This gives pointwise inequality between measures
\begin{equation}\label{major}
\mathcal{C}_1+\cdots+\mathcal{C}_s\leq \mathcal{B}_1+\cdots+\mathcal{B}_s,\quad s\geq 1.
\end{equation}
Applying the Markov inequality,  the total variation distance is estimated in terms of  the means as 
\begin{align*}\label{TV-pro}
&\hspace{-1cm}d_{\rm TV}(
(\mathcal{C}_s((n_0,n_1]), 1\leq s\leq r), (\mathcal{B}_s((n_0,n_1]), 1\leq s\leq r))\\
&\leq {\mathbb P} [\mathcal{C}_s((n_0,n_1]), 1\leq s\leq r)\neq (\mathcal{B}_s((n_0,n_1]), 1\leq s\leq r)) ]\\
&= {\mathbb P} \left[\sum_{s=1}^{r_1}\mathcal{C}_s((n_0,n_1])\neq \sum_{s=1}^{r_1}\mathcal{B}_s((n_0,n_1])\;\;{\rm for}\;\; {\rm some}\;\; r_1\leq r \right]\\
&\leq \sum_{r_1=1}^{r} {\mathbb P} \left[ \sum_{s=1}^{r_1}    \left( \mathcal{B}_s((n_0,n_1])-\mathcal{C}_s((n_0,n_1]) \right) \geq 1\right]\\
&\leq \sum_{r_1=1}^{r} {\mathbb E} \left[ \sum_{s=1}^{r_1}\{ \mathcal{B}_s((n_0,n_1])-\mathcal{C}_s((n_0,n_1]) \} \right].
\end{align*}
From (\ref{mean1}), (\ref{mean2}), the bound approaches zero as $n_0\to\infty$ and $n_1=n_0+\lfloor \theta n_0/\Phi(n_0)\rfloor$, provided that $\Phi$ satisfies the condition of regular or slow variation with $\alpha\in [0,1)$. For the sequel we assume this holds, and focus on approximating $({\mathcal{C}_s}((n_0,n_1]), ~1\leq s\leq r)$.

\vspace{2mm}
\noindent
{\bf Step 2:} Multinomial approximation to the conditional distribution of $(\mathcal{C}_s((n_0,n_1]),~1\leq s\leq r)$ given $\{\mathcal{K}_{n_1,s}=k_s,~1\leq s\leq r\}$. 

\vspace{2mm}
The processes $\mathcal{C}_1,\ldots, \mathcal{C}_r$ do not have common points and have the following structure. Conditionally on 
$\mathcal{K}_{n_1,s}=k_s,~1\leq s\leq r,$
they altogether  can be represented by sampling without replacement from an urn with $s k_s$ balls of colour $s$ (that occupy $s$-ton boxes at time $n_1$)
and $n_1-\sum_{s=1}^r s k_s$ uncoloured balls. In particular, the conditional distribution of the random vector $(\mathcal{C}_s((n_0,n_1]),~1\leq s\leq r)$ (complemented with the uncoloured component) given $\{\mathcal{K}_{n_1,s}=k_s,~1\leq s\leq r\}$ is a multivariate hypergeometric distribution with parameters
$$
n_1-n_0;\, \left(k_1, 2k_2,\ldots, r k_r, n_1-\sum_{s=1}^r s k_s\right).
$$
By the Diaconis-Freedman bound, see Theorem 4 in~\cite{Diaconis}, the total variation distance between this multivariate hypergeometric distribution and its multinomial counterpart with parameters
\begin{equation}\label{HG-par}
n_1-n_0;\, \left(\frac{k_1}{n_1},\frac{2k_2}{n_2},\ldots,\frac{rk_r}{n_1},\frac{n_1-\sum_{s=1}^r s k_s}{n_1}\right)
\end{equation}
is at most $2(r+1)(n_1-n_0)/n_1$. Remarkably, the bound does not depend on $k_1,\ldots,k_r$.
{According to Theorem 1 in~\cite{McDonald} the total variation distance between the first $r$ components of this multinomial distribution and the $r$-variate Poisson distribution 
$${\rm Poiss}\left(\frac{(n_1-n_0)}{n_1}\,s k_s,~1\leq s\leq r\right)$$
is bounded by $2(n_1-n_0)\left(\sum_{s=1}^{r}\frac{s k_s}{n_1}\right)^2$ which goes to zero.}

\vspace{2mm}
\noindent
{\bf Step 3:} Mixed multivariate Poisson approximation of $(\mathcal{C}_s((n_0,n_1]),~1\leq s\leq r)$.

\vspace{2mm}
To eliminate conditioning  we use the elementary consequence of the fact that the total variation distance derives from a norm. This fact is implicit in~\cite{Diaconis}, see the proof of Theorem 3 therein.
\begin{lemma}
For  two families of probability measures $(F_\alpha), (G_\alpha)$ the convex mixtures satisfy
$$d_{\rm TV}\left(\sum_\alpha a_\alpha F_\alpha, \sum_\alpha a_\alpha G_\alpha\right)\leq \sum_\alpha  a_\alpha~ d_{\rm TV}( F_\alpha, G_\alpha).$$
\end{lemma}
\noindent
Together with the estimates from Steps 2 and 3 
this allows us to assess the approximation 
by a mixed multivariate Poisson distribution with random parameters:
\begin{multline}
d_{\rm TV}\left( (\mathcal{C}_s((n_0,n_1]),~1\leq s\leq r),~{\rm Poiss}\left(\frac{(n_1-n_0)}{n_1} ~s \mathcal{K}_{n_1,s}, ~1\leq s\leq r\right) \right)\\
\leq \frac{(r+1)(n_1-n_0)}{n_1}+\frac{2(n_1-n_0)}{n_1^2}{\mathbb E}\left[ \left(\sum_{s=1}^r s \mathcal{K}_{n_1,s}\right)^2\right].\label{mixedP}
\end{multline}
To proceed we need the following lemma. For the Poissonised scheme it is  quite easy too see that ${\rm Var}[K_r(t)]<{\mathbb E}[K_r(t)]=\Phi_r(t)$; see p.~370 in~\cite{GBarbour}. We are not aware of discrete time analogue of this inequality, and will use instead a rougher estimate.
\begin{lemma}\label{lem:variance} For $r\geq 1$ and $n\geq 1$,
$$
{\rm Var}\left[\sum_{s=1}^r \mathcal{K}_{n,s}\right]\leq {\mathbb E}\left[\sum_{s=1}^r \mathcal{K}_{n,s}\right]\quad\text{and}\quad {\rm Var}\left[ \mathcal{K}_{n,r}\right]\leq 4 \,{\mathbb E}\left[\sum_{s=1}^r \mathcal{K}_{n,s}\right].
$$
\end{lemma} 
\proof The array $(\mathcal{P}_j(n),~j\geq 1)$ has multinomial distribution, hence it is negatively associated; see Section 3 in~\cite{Joag-Dev}. The indicators $\mathbbm{1}_{\{\mathcal{P}_j(n)\leq r\}}$ are nonincreasing functions of the array hence they are pairwise negatively correlated. Thus,
$$
{\rm Var}\left[\sum_{s=1}^r \mathcal{K}_{n,s}\right]={\rm Var}\left[\sum_{j\geq 1} \mathbbm{1}_{\{\mathcal{P}_j(n)\leq r\}}\right]\leq \sum_{j\geq 1} {\rm Var}[\mathbbm{1}_{\{\mathcal{P}_j(n)\leq r\}}]\leq  \sum_{j\geq 1} {\mathbb E}[\mathbbm{1}_{\{\mathcal{P}_j(n)\leq r\}}]={\mathbb E}\left[\sum_{s=1}^r \mathcal{K}_{n,s}\right].
$$
By the virtue of
$$
\mathcal{K}_{n,r}=\sum_{s=1}^r \mathcal{K}_{n,s}- \sum_{s=1}^{r-1} \mathcal{K}_{n,s}
$$
the second inequality follows from the first.
\endpf

By the first inequality in Lemma~\ref{lem:variance} the second summand on the right-hand side of~\eqref{mixedP} is bounded by
$$
\frac{2r^2 (n_1-n_0)}{n_1^2}{\mathbb E}\left[ \left(\sum_{s=1}^r \mathcal{K}_{n_1,s}\right)^2\right]\leq \frac{2r^2(n_1-n_0)}{n_1^2}\left({\mathbb E}\left[\sum_{s=1}^r \mathcal{K}_{n_1,s}\right]+\left({\mathbb E}\left[\sum_{s=1}^r \mathcal{K}_{n_1,s}\right]\right)^2\right).
$$
Since ${\mathbb E}\left[\sum_{s=1}^r \mathcal{K}_{n_1,s}\right]=\sum_{s=1}^r \left(\Phi_s(n_1)+o(1)\right)$ and $n_1=n_0+\lfloor\theta f(n_0)\rfloor$, this tends to zero in case $\alpha\in [0,1)$ of Theorem \ref{T2}.

\vspace{2mm}
\noindent
{\bf Step 4:} Multivariate Poisson approximation.

\vspace{2mm}
Finally, we wish to replace the random parameters in (\ref{mixedP}) by their mean values. Applying Lemma~\ref{lem:variance} and Theorem 10.C from~\cite{Barbour}
\begin{align*}
&\hspace{-1cm}d_{\rm TV}\left({\rm Poiss}\left(\frac{(n_1-n_0)s}{n_1}  \mathcal{K}_{n_1,s}, ~1\leq s \leq r\right), {\rm Poiss}\left(\frac{(n_1-n_0)s}{n_1}  {\mathbb E}[\mathcal{K}_{n_1,s}], ~1\leq s\leq r\right) \right)\\
&\leq{\mathbb E}\left[\left( \sum_{s=1}^r \frac{(n_1-n_0)s}{n_1}\left|\mathcal{K}_{n_1,s}- {\mathbb E}[\mathcal{K}_{n_1,s}] \right|\right)^2\right]\leq
 \frac{r(n_1-n_0)^2}{n_1^2}\sum_{s=1}^r s^2 \,{\rm Var}[\mathcal{K}_{n_1,s}]\\
&\leq \frac{4r^4(n_1-n_0)^2}{n_1^2}\sum_{s=1}^r {\mathbb E}[\mathcal{K}_{n_1,s}]=
\frac{4r^4(n_1-n_0)^2}{n_1^2}\sum_{s=1}^r \left(\Phi_s(n_1)+o(1)\right).
\end{align*}
For $n_1=n_0+\lfloor\theta f(n_0)\rfloor$ this bound approaches $0$ in all cases of Theorem \ref{T2}. In particular, in the proper regular variation  case (i.e. with index $0<\alpha<1$) with the scaling function 
 $f(n_0)=\Phi(n_0)$ the bound is of the order
$O(1/\Phi(n_0))$. The same bound is valid   if we approximate by ${\rm Poiss}(\theta\rho_s,~1\leq s\leq r)$, that is with the interpolated
rate instead of the natural mean.

From these total variation bounds the analogue of Theorem \ref{T2} for  $(\widehat{\mathcal B}_r(\theta), ~1\leq s\leq r)$ follows in case $\alpha\in[0,1)$. In particular, in the proper regular variation case, the vector 
converges in distribution to $({\rm Poiss}(\theta\rho_s), 1\leq s\leq r)$.

\section{Random frequencies}

Finally, we  sketch a mixed Poisson approximation for records in occupancy schemes where the frequencies $(p_j)$ are random. For simplicity of exposition we shall consider only the classic `pure power laws', see Section 10 in~\cite{GHP}. More precisely, assume that for some $\alpha\in (0,1)$ 
\begin{equation}\label{eq:random_freq_nu}
\nu[x,1]~\sim~D x^{-\alpha},\quad x\to 0+\quad\text{a.s.},
\end{equation}
where $D$ is a strictly positive random variable sometimes called the $\alpha$-diversity \cite{CSP}.
Intensely studied examples of $(p_j)$ leading to (\ref{eq:random_freq_nu}) are the two-parameter Poisson-Di\-ri\-chlet frequencies \cite{CSP} and their generalisations \cite{Gibbs, GPY, Spano}.
By Proposition 23 in~\cite{GHP} the relation 
(\ref{eq:random_freq_nu})
is equivalent to
\begin{equation}\label{eq:random_freq_p}
p_j~\sim~ D^{1/\alpha}j^{-1/\alpha},\quad j\to\infty\quad\text{a.s.}
\end{equation}
By Karamata's theorem, see Theorem 1.6.1 in~\cite{Goldie}, either of these relations implies
\begin{equation}\label{eq:random_freq_laplace}
\Phi(n_0)~\sim~\Gamma(1-\alpha) D n_0^{\alpha},\quad n_0\to \infty\quad\text{a.s.}
\end{equation}
By Theorem~\ref{T4}, the processes
$$
\left(\widehat{\mathcal{B}}_r\left(\left(n_0,n_0+\theta\frac{n_0}{\Phi(n_0)}\right]\right),\,\theta\geq 0\right),\quad r\geq 1,
$$
given $(p_j)$, converge jointly, as $n_0\to\infty$, to independent homogeneous Poisson processes with rates $\rho_r=\frac{\alpha \Gamma(r-\alpha)}{(r-1)!\, \Gamma(1-\alpha)}$. Combining this with~\eqref{eq:random_freq_laplace} we arrive at
\begin{prop}
Assume either of equivalent conditions~\eqref{eq:random_freq_nu} or~\eqref{eq:random_freq_p}. Then, as $n_0\to\infty$,
$$
\left(\widehat{\mathcal{B}}_r\left(\left(n_0,n_0+n_0^{1-\alpha}\theta\right]\right),\,\theta\geq 0\right),
$$
converge jointly in distribution as $n_0\to\infty$, 
to mixed
Poisson processes with random rates $\frac{D\alpha \Gamma(r-\alpha)}{(r-1)!}$. 
Conditionally on  $D$, the limit processes are independent Poisson.
\end{prop}

\section*{Appendix}

\begin{lemma}\label{L-app} For $\mu$ a measure on the halfline with $0<\mu([0,\infty))\leq 1$, the 
 Laplace transform
$$L(t):=\int_0^\infty e^{-tx}\mu({\rm d}x),\quad t\geq 0,$$
satisfies  for $0\leq\gamma\leq 1$
$$L(\gamma t)\leq      \,L^\gamma (t).$$
\end{lemma}

\proof  Let $\xi$ be a nonnegative random variable with Laplace transform $L(t)/L(0)$. Using Jensen's inequality,
$$\frac{L(\gamma t)}{L(0)}={\mathbb E}[e^{-\gamma t\xi}]\leq \left({\mathbb E}[e^{-t\xi}]\right)^\gamma=\left( \frac{L(t)}{L(0)}\right)^\gamma.$$
The assertion follows from this by noting that $L(0)=\mu([0,\infty))\leq 1$.
\endpf

\vspace{2mm}
\noindent
{\bf Remark on exchangeable partitions.} 
Nacu \cite{Nacu} proved that the distribution of  the point process of records $\mathcal{B}_1$ restricted to $\{1,\ldots,n\}$ uniquely determines the 
distribution of counts $({\mathcal K}_{n,1},\ldots,{\mathcal K}_{n,n})$ that encode a partition of integer $n$ induced by the allocation of $n$ balls.
Letting $n$ vary, by the virtue of Kingman's theory of exchangeable partitions \cite{CSP} this fact implies that the distribution of probabilities $(p_j)$ (arranged in decreasing order) is uniquely determined by the distribution of $\mathcal{B}_1$ seen as 
a point process on ${\mathbb N}$.

Recall that ${\mathcal C}_r(n_1)=r{\mathcal K}_{n_1,r}$. Changing notation $n_1$ to $n$
 we have from (\ref{major}) 
$$
\mathcal{K}_{n,1}+ 2 \mathcal{K}_{n,2}+  \cdots+r\mathcal{K}_{n,r}\leq \mathcal{B}_1(n)+\mathcal{B}_2(n)+\cdots+\mathcal{B}_r(n),\quad    1\leq r \leq n.
$$
Now, using this inequality, arguments similar  to \cite{Nacu} allow one to show that the distribution of
$(\mathcal{B}_1(n),\ldots,\mathcal{B}_n(n))$ uniquely determines the distribution of the partition of $n$.
Letting $n$ vary, we conclude that the distributions of  vectors $(\mathcal{B}_1(n),\ldots,\mathcal{B}_n(n)), n\geq 1,$ offer yet another way to describe the law of the exchangeable partition  of ${\mathbb N}$ induced by the allocation of infinitely many balls.

\section*{Acknowledgments}

This work has been accomplished during AM's visit to Queen Mary University of London as Leverhulme Visiting Professor in July-December 2023. AM gratefully acknowledges financial support from the Leverhulme Trust.


\begin{thebibliography}{99}


\bibitem{GBarbour} A.~Barbour and A.~Gnedin (2009).  Small counts in the infinite occupancy scheme, {\it Electron. J. Probab.} {\bf 14}, Paper no. 13, 365--384.

\bibitem{Barbour} A.~Barbour, L.~Holst and S.~Jansen (1992). {\it  Poisson approximation}, Clarendon Press, Oxford.

\bibitem{Ben} A.~Ben-Hamou, S.~Boucheron and M.~I.~Ohannessian (2017).
Concentration inequalities in the infinite urn scheme for occupancy counts and the missing mass, with applications, {\it Bernoulli} {\bf 23}, 249--287.


\bibitem{Goldie} N.~H.~Bingham, C.~M.~Goldie and J.~L.~Teugels (1987). {\it Regular variation}, Cambridge University Press.


\bibitem{Bunge} J.~Bunge and M.~Fitzpatrick (1993). Estimating the number of species: A Review, 
{\it J. Amer. Statist. Assoc.} {\bf  88}, 364--373.


\bibitem{Serfozo} F.~B{\"o}ker and R.~Serfozo (1983). Ordered thinnings of point processes and random measures, {\it Stoch. Proc.  Appl.} {\bf  15}, 113--132.

\bibitem{Klesov} V.~V.~Buldygin, K.-H.~Indlekofer, O.~I.~Klesov and I.~Steinebach (2018). {\it Pseudo-regularly varying functions and generalized renewal processes}, Springer.

\bibitem{Moderate} D.~Buraczewski, B.~Dovgay and A.~Marynych (2021). Moderate parts in regenerative compositions: the case of regular variation, {\it J. Math. Analysis Appl.} {\bf 497}, 124894.

\bibitem{Variance} L.~Bogachev, A.~Gnedin and Y.~Yakubovich (2008). On the variance of the number of occupied boxes, {\it Adv. Appl. Math.} {\bf 40}, 401--432.

\bibitem{CZ} M.~Chebunin and S.~Zuyev (2022). Functional central limit theorems for occupancies and missing mass process in infinite urn models, {\it J. Theor. Probab.} {\bf 35}, 1--19.

\bibitem{Daley} D.~J.~Daley and D.~Vere-Jones (2008). {\it An introduction to the theory of point processes. Volume II: general theory and structure}, Springer.

\bibitem{Diaconis} P.~Diaconis and D.~Freedman (1980). Finite exchangeable sequences, {\it Ann. Probab.} {\bf 8}, 745--764.

\bibitem{DW} O.~Durieu and Y.~Wang (2016). From infinite urn schemes to decompositions of self-similar Gaussian processes, {\it Electron. J. Probab.} {\bf 21}, Paper no. 43, 1--23.

\bibitem{GHP} A.~Gnedin, B.~Hansen and J.~Pitman (2007).  Notes on the occupancy problem with infinitely many boxes: general asymptotics and power laws, {\it Probab. Surv. } {\bf 4}, 146--171.


\bibitem{Gibbs} A.~Gnedin and J.~Pitman (2006). Exchangeable Gibbs partitions and Stirling triangles,
{\it J. Math. Sci.} {\bf  138}, 5674--5685.

\bibitem{GPY} A.~Gnedin, J.~Pitman and M. Yor (2006).
Asymptotic laws for compositions derived from transformed subordinators, {\it Ann. Probab.} {\bf 34}, 468--492.


\bibitem{GS} A.~Gnedin and D.~Stark (2023). Random permutations and queues, {\it Adv. Appl. Math.} {\bf 149}, 102549.

\bibitem{Spano} R.~C.~Griffiths and D.~Span{\'o} (2007). Record indices and age-ordered frequencies in exchangeable Gibbs partitions, {\it Electron. J. Probab.} {\bf 12}, 1101--1130.

\bibitem{Hardy} G.~H.~Hardy (1949). {\it Divergent series}, Clarendon Press, Oxford.

\bibitem{Hwang} H.-K.~Hwang and S.~Janson (2008). Local limit theorems for finite and infinite urn models, {\it Ann. Probab.} {\bf 36}, 992--1022.

\bibitem{Ilienko} A. Ilienko (2019). Convergence of point processes associated with coupon collector's and Dixie cup problems, {\it Electron. Commun. Probab.} {\bf 24}, Paper no. 51, 1--9.


\bibitem{Joag-Dev} K.~Joag-Dev and F.~Proschan (1983). Negative association of random variables with applications, {\it Ann. Statist.} {\bf 11}, 286--295.

\bibitem{Kallenberg} O.~Kallenberg (2017). {\it Random measures: theory and applications}, Springer.

\bibitem{Karlin} S.~Karlin (1967). Central limit theorems for certain infinite urn schemes, {\it J. Math. Mech.} {\bf 17}, 373--401.


\bibitem{Last} G.~Last and R.~Szekli (2019). On negative association of some finite point processes on general state spaces, {\it J. Appl. Probab.}, {\bf 56}, 139--152.

\bibitem{McDonald} D.~R.~McDonald (1980). On the Poisson approximation to the multinomial distribution, {\it Canadian J.  Statist.} {\bf 8}, 115--118.

\bibitem{Nacu} S.~Nacu (2006). Increments of reandom partitions, {\it Comb. Prob. Comp.} {\bf 15},  589--595.

\bibitem{CSP} J.~Pitman, {\it Combinatorial stochastic processes}, Springer Lecture Notes in Mathematics {\bf 1875}.
\bibitem{Resnick} S.~Resnick (2008). {\it Extreme values, regular variation and point processes}, Springer.

\bibitem{Roos} B.~Roos (1999). On the rate of multivariate Poisson convergence,
{\it J. Multivariate Anal.} {\bf 69}, 120--134.

\bibitem{Ruzankin} P.~S.~Ruzankin (2004). On the rate of Poisson process approximation to a Bernoulli process, {\it J. Appl. Probab.} {\bf 41}, 271--276.

\bibitem{Schilling} R.~L.~Schilling, R.~Song and Z.~Vondra\v{c}ek (2012). {\it Bernstein functions}, De Gruyter.


\end{thebibliography}
\end{document}